\documentclass{amsart}
\usepackage{amsfonts,amsmath,amsthm,amscd,amssymb,latexsym}
\usepackage{epsfig}

\begin{document}

\title{On Sushchansky $p$-groups}
\author{Ievgen V. Bondarenko}
\address[I.~Bondarenko]{Kyiv Taras Shevchenko University, Ukraine and Department of Mathematics, Texas A\&M University, MS-3368, College Station, TX 77845-3368, USA}
\email{ibond@math.tamu.edu}

\author{Dmytro M. Savchuk}
\thanks{Both authors were partially supported by NSF grants DMS-0308985 and
DMS-0456185}
\address[D.~Savchuk]{Department of Mathematics, Texas A\&M University, MS-3368, College Station, TX 77845-3368, USA}
\email{savchuk\@math.tamu.edu}
\urladdr{www.math.tamu.edu/\~{}savchuk}

\subjclass[2000]{20F69, 20F10, 20E08} \keywords{Burnside groups,
growth of groups, automata groups, branch groups}

%\dedicated{Dedicated to V.I.~Sushchansky on the occasion of his 60th
%birthday}
%\communicated{Communicated person1}

%\research% Type of the paper. \survey command may be used also.

%\date{06.10.2006}

%Change theorem environments according to your needs...
\newtheorem{thm}{Theorem}
\newtheorem{prop}[thm]{Proposition}
\newtheorem{cor}[thm]{Corollary}
\newtheorem{lemma}[thm]{Lemma}
\theoremstyle{definition}
\newtheorem{defin}{Definition}
\newtheorem{exam}{Example}

\sloppy \frenchspacing

\newcommand{\xs}{X^*}
\newcommand{\A}{\mathcal{A}}
\newcommand{\Aut}{\mathop{\rm Aut}\nolimits}
\newcommand{\Sym}{\mathop{\rm Sym}\nolimits}
\newcommand{\Z}{\mathbb Z}
\newcommand{\G}{\mathsf{G}}
\newcommand{\Fp}{\mathbb{F}_p}
\newcommand{\St}{\mathop{\rm St}\nolimits}
\newcommand{\Orb}{\mathop{\rm Orb}\nolimits}

\maketitle

\begin{center}
\textit{Dedicated to V.I.~Sushchansky on the occasion of his 60th
birthday}
\end{center}

\begin{abstract}
We study Sushchansky $p$-groups introduced in
\cite{sushch:burnside}. We recall the original definition and
translate it into the language of automata groups. The original
actions of Sushchansky groups on $p$-ary tree are not
level-transitive and we describe their orbit trees. This allows us
to simplify the definition and prove that these groups admit
faithful level-transitive actions on the same tree. Certain branch
structures in their self-similar closures are established. We
provide the connection with, so-called, $\G$
groups~\cite{bar_gs:branch} that shows that all Sushchansky groups
have intermediate growth and allows to obtain an upper bound on
their period growth functions.
\end{abstract}

%\subjclass{2001}{20F69, 20F10, 20E08}

%\keywords{Burnside groups, growth of groups, automata groups, branch
%groups}

\section*{Introduction}

Sushchansky $p$-groups were introduced in~\cite{sushch:burnside} as
one of the pioneering examples of finitely generated infinite
torsion groups, providing counter-examples to the General Burnside
problem. Initially, this problem was solved by E.S.~Golod
in~\cite{golod:p-groups} using the Golod-Shafarevich theorem.
Simpler and easier to handle counter-examples were constructed by
S.V.~Aleshin in~\cite{aleshin:burnside} by means of automata. The
use of automata groups to resolve Burnside's problem was earlier
suggested by V.M.~Glushkov in~\cite{glushkov:automata}. But only
after the results of R.I.~Grigorchuk
from~\cite{grigorch:burnside,grigorch:milnor} automata groups became
the subject of deeper investigation. It happened that this class
contains groups with many extraordinary properties, like infinite
torsion groups, groups of intermediate growth, groups of finite
width, just-infinite groups, etc.

V.I.~Sushchansky used a different language, namely the language of
tableaux, introduced by L.~Kaluzhnin to study properties of iterated
wreath products~\cite{kalou:la_structure}. For each prime $p>2$,
V.I.~Sushchansky constructed a finite family of infinite $p$-groups
generated by two tableaux. Each such a tableau naturally defines an
automorphism of a rooted tree and, as was already noticed in
\cite{gns00:automata}, can be represented by a finite initial
automaton. We describe these automata and study Sushchansky groups
and their actions on rooted trees by means of this well-developed
language.

The structure of the paper is as follows. In
Section~\ref{sec_original} we recall the original definition of
Sushchansky groups. In Section~\ref{sec_automata} we describe the
corresponding automata. The associated action on a rooted tree is
not level-transitive and in Section~\ref{sec_trees} we describe its
orbit tree and show that there exists a faithful level-transitive
action given by finite initial automata. The self-similar closure is
studied in Section~\ref{sec_closure}. The main results are presented
in Section~\ref{sec_growth}. It was pointed out
in~\cite{grigorch:degrees85} that all Sushchansky $p$-groups have
intermediate growth, but only the main idea of the proof was given.
Here we provide a complete proof of this fact together with new
estimates on the growth function, thus contributing to the Milnor
question \cite{milnor:problem}, which was solved
in~\cite{grigorch:milnor} by R.I.~Grigorchuk. Also we give an upper
bound on the period growth function. The main idea is to use $\G$
groups of intermediate growth introduced
in~\cite{bartholdi_s:growth} (see also~\cite{bar_gs:branch}). For
each Sushchansky $p$-group we construct a $\G$ group of intermediate
growth and prove that their growth functions are equivalent.

The authors wish to thank Zoran \v Suni\'c for fruitful discussions
and important comments, which enhanced the paper.

\section{Original definition via tableaux}\label{sec_original}
Let $X=\{0,1,\ldots, p-1\}$ be a finite alphabet for some prime $p$.
We identify $X$ with the finite field $\Fp$.

The set $X^{*}$ of all finite words over $X$ has a natural structure
of a rooted $p$-ary tree. Every automorphisms $g\in\Aut X^{*}$ of
this tree induces an automorphism $g|_v$ of the subtree $vX^{*}$ by
the rule $g|_v(w)=u$ if and only if $g(vw)=g(v)u$. This automorphism
is called \textit{the restriction of $g$ on word $v$} (in some
papers the word \textit{section} or \textit{state} is used).

The Sylow $p$-subgroup $P_{\infty}$ of the profinite group $\Aut
X^*$ is equal to the infinite wreath product of cyclic groups of
order $p$, i.e. $P_{\infty}=\wr_{i\geq 1} C_p^{(i)}$. Using this
description one can construct special ``tableau'' representation of
$P_{\infty}$. The ``tableau'' representation was initially
introduced by L.~Kaluzhnin for Sylow $p$-subgroups of symmetric
groups of order $p^m$ in~\cite{kalou:la_structure}.

The group $P_{\infty}$ is isomorphic to the group of triangular
tableaux of the form:
\[
u=[a_1,a_2(x_1), a_3(x_1,x_2),\ldots],
\]
where $a_1\in\Fp$,
$a_{i+1}(x_1,\ldots,x_i)\in\Fp[x_1,\ldots,x_i]/\langle
x_1^p-x_1,\ldots, x_i^p-x_i\rangle$. The multiplication of tableaux
is given by the formula:
\[
[a_1, a_2(x_1), a_3(x_1,x_2),\ldots ]\cdot[b_1, b_2(x_1),
b_3(x_1,x_2), \ldots]=
\]
\[
 =[a_1+b_1,a_2(x_1)+b_2(x_1+a_1), a_3(x_1,x_2)+b_3(x_1+a_1,x_2+a_2(x_1)), \ldots].
\]

\noindent The action of the tableau $u$ on the tree $X^{*}$ is given
by:
\begin{equation}\label{eqn tableau action}
u(x_1x_2\ldots x_n)=y_1y_2\ldots y_n,
\end{equation}
where $y_1=x_1+a_1$, $y_2=x_2+a_2(x_1),\ldots,
y_n=x_n+a_n(x_1,\ldots,x_{n-1})$, where all calculations are made by
identifying $X$ with the field $\Fp$.

For the duration of the rest of the paper we fix a prime $p>2$.

Fix some order $\lambda=\{(\alpha_i,\beta_i), i=1,\ldots, p^2\}$ on
the set of pairs $\{(\alpha,\beta) | \alpha,\beta\in \Fp\}$. For
$j>p^2$ we define $(\alpha_j,\beta_j)=(\alpha_i,\beta_i)$ where
$i\equiv j \mod p^2$. Define two tableaux
\[
A=[1,x_1,0,0,\ldots], \ \
B_{\lambda}=[0,0,b_3(x_1,x_2),b_4(x_1,x_2,x_3),\ldots],
\]
where the coordinates of $B_{\lambda}$ are defined by its values in
the following way:
\begin{enumerate}
    \item[a)] $b_3(2,1)=1$;
    \item[b)] $b_i(0,0,\ldots,0,1)=1$ if $\beta_i\neq0$;
    \item[c)] $b_i(1,0,\ldots,0,1)=-\frac{\alpha_i}{\beta_i}$ if $\beta_i\neq0$ and
    $b_i(1,0,\ldots,0,1)=1$ if $\beta_i=0$;
    \item[d)] all the other values are zeroes.
\end{enumerate}

The group $G_{\lambda}=\langle A, B_{\lambda}\rangle$ is called
\textit{the Sushchansky group of type $\lambda$}. The following
theorem is proven in \cite{sushch:burnside}.

\begin{thm} \label{thm_susch}
$G_{\lambda}$ is infinite periodic $p$-group for any type $\lambda$.
\end{thm}

\section{Automata approach}
\label{sec_automata}

Another language dealing with groups acting on rooted trees is the
language of automata groups. For a definitions we refer to the
survey paper~\cite{gns00:automata}. Many groups related to Burnside
and Milnor Problems happen to be in the class of groups generated by
finite automata. The Sushchansky groups are not an exception and we
describe the structure of the corresponding automata in this
section.

The action of every automorphism $g$ of the rooted tree $X^*$ can be
encoded by an initial automaton whose states are the restrictions of
$g$ on the finite words over $X$. In the case when this set is
finite we call $g$ a \emph{finite-state} automorphism. The action of
such an automorphism is encoded by a finite automaton.

It is known that (see~\cite{gns00:automata}) that $\Aut X^*\cong\Aut
X^*\wr\Sym(X)$, which gives a convenient way to represent every
automorphism in the following form:
\[g=(g|_0,g|_1,\ldots,g|_{p-1})\pi_g,\]
where $g|_0,g|_1,\ldots,g|_{p-1}$ are the restrictions of $g$ on the
letters of $X$ and $\pi_g$ is the permutation of $X$ induced by $g$.

The multiplication of automorphisms written in this way is performed
as follows. If $h=(h|_0,h|_1,\ldots,h|_{p-1})\pi_h$ then
\[gh=(g|_0h|_{\pi_g(0)},\ldots,g|_{p-1}h|_{\pi_g(p-1)})\pi_g\pi_h.\]

Now we proceed with an explicit construction of automata associated
to Sushchansky groups. Let $\sigma=(0,1,\ldots,p-1)$ be a cyclic
permutation of $X$. With a slight abuse of notation, depending on
the context, $\sigma$ will also denote the automorphism of $X^*$ of
the form $(1,1,\ldots,1)\sigma$.

Given the order $\lambda=\{(\alpha_i,\beta_i)\}$ define words
$u,v\in X^{p^2}$ in the following way:
\[
u_i=\left\{%
\begin{array}{ll}
    0, & \hbox{ if } \beta_i=0; \\
    1, & \hbox{ if } \beta_i\neq 0. \\
\end{array}%
\right. \qquad \qquad
v_i=\left\{%
\begin{array}{ll}
    1, & \hbox{ if } \beta_i=0; \\
    -\frac{\alpha_i}{\beta_i}, & \hbox{ if } \beta_i\neq 0. \\
\end{array}%
\right.
\]
The words $u$ and $v$ encode the actions of $B_{\lambda}$ on the
words $00\ldots 01*$ and $10\ldots 01*$, respectively. Using the
words $u$ and $v$ we can construct automorphisms
$q_1,\ldots,q_{p^2}, r_1,\ldots, r_{p^2}$ of the tree $X^{*}$ by the
following recurrent formulas:
\begin{equation}\label{eqn def aut q r}
q_i=(q_{i+1},\sigma^{u_i},1,\ldots,1), \qquad
r_i=(r_{i+1},\sigma^{v_i},1,\ldots ,1),
\end{equation}
for $i=1,\ldots,p^2$, where the indices are considered modulo $p^2$,
i.e. $i=i+np^2$ for any $n$.

Formula (\ref{eqn tableau action}) implies that $q_i$ and $r_i$ are
precisely the restrictions of $B_{\lambda}$ on the words
$00(0)^{i-1+np^2}$ and $10(0)^{i-1+np^2}$, respectively, for any
$n\geq 0$.

\begin{figure}[h]
\begin{center}
\epsfig{file=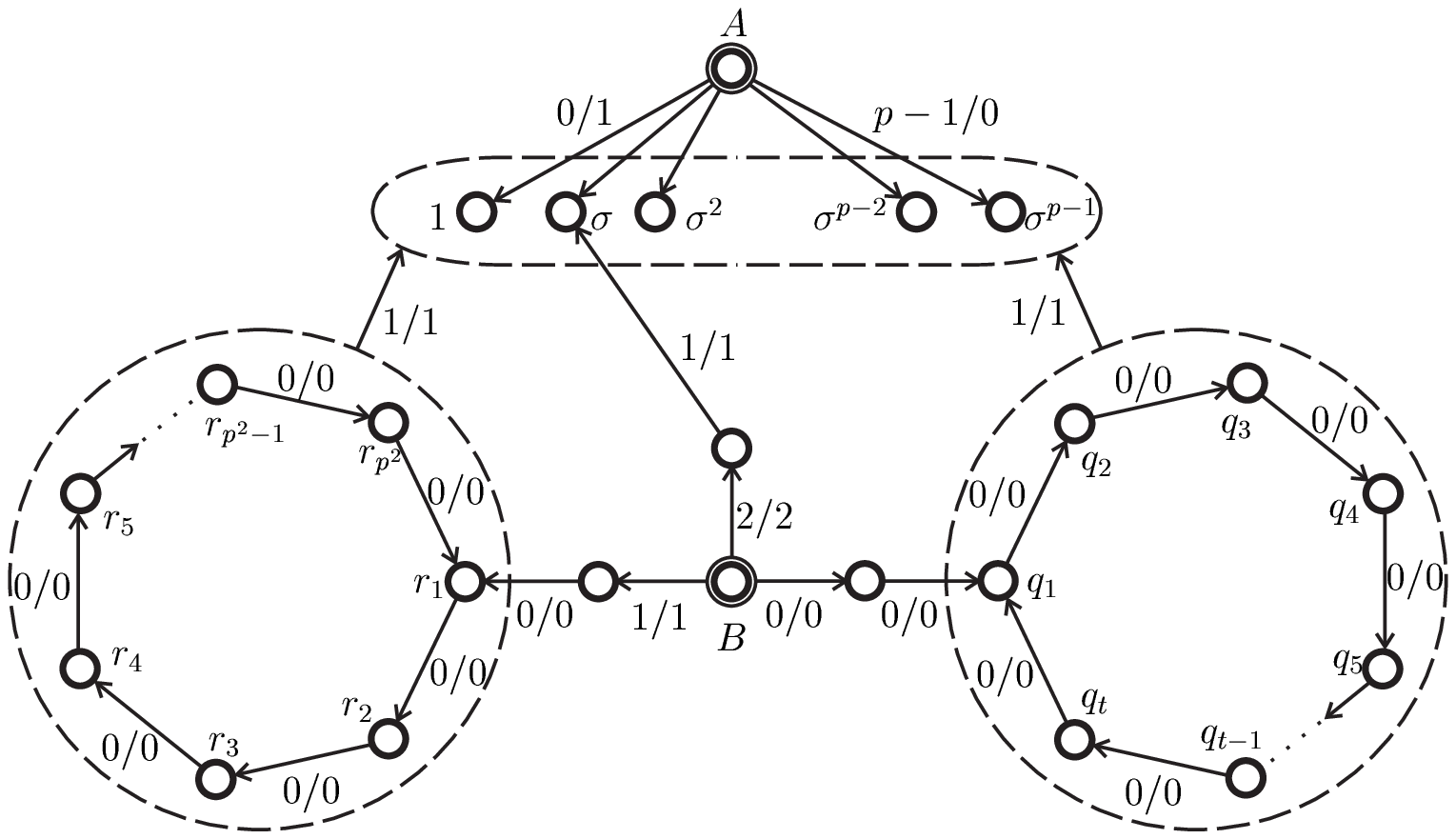,width=340pt}
\end{center}
\caption{The Structure of Sushchansky automata\label{aut_general}}
\end{figure}

\begin{figure}[h]
\begin{center}
\epsfig{file=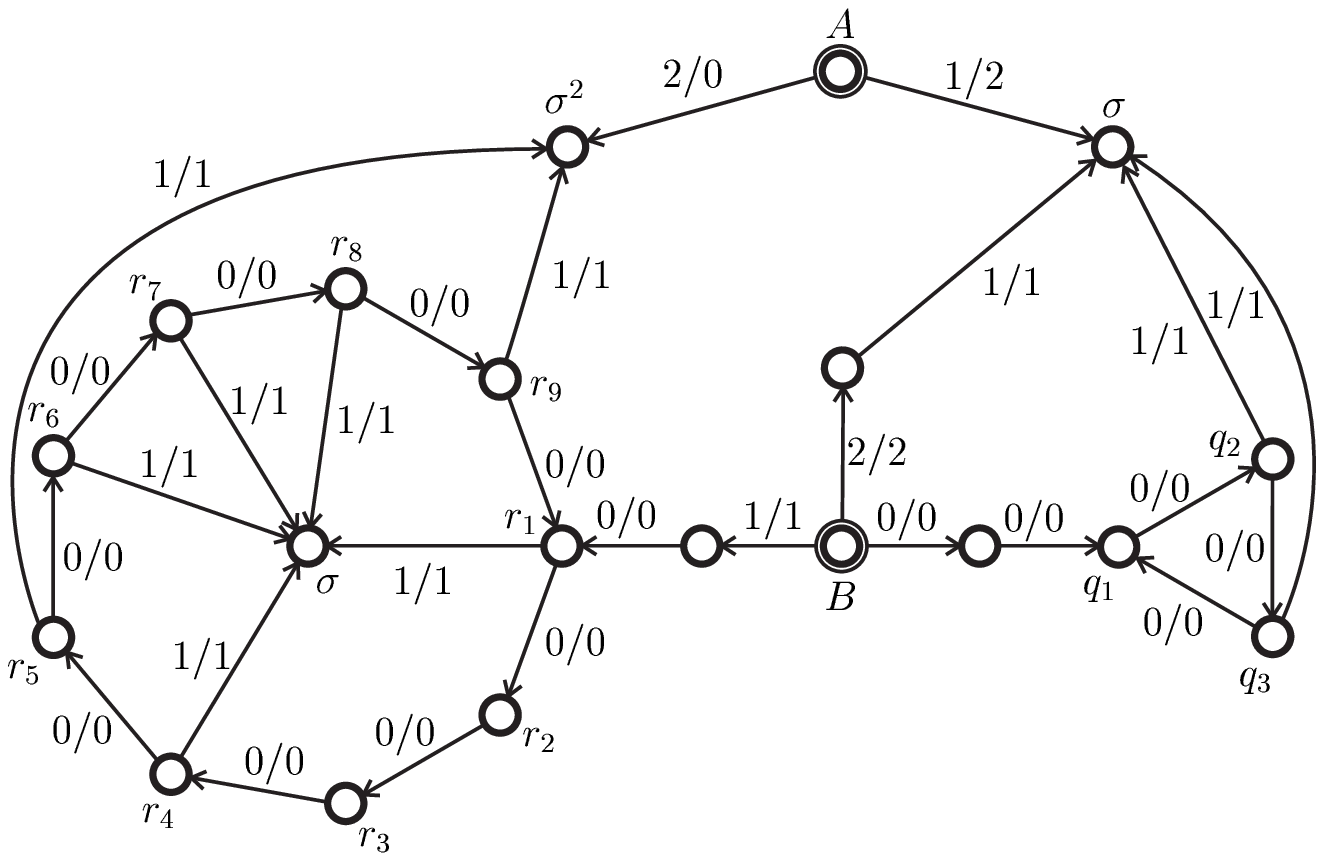,width=340pt}
\end{center}
\caption{Sushchansky automaton for $p=3$ corresponding to the
lexicographic order\label{aut_particular}}
\end{figure}

The action of the tableau $A$ is given by:
\[
A=(1,\sigma,\sigma^2,\ldots,\sigma^{p-1})\sigma;
\]
while $B_{\lambda}$ acts trivially on the second level and the
action on the rest is given by the restrictions:
\[
B_{\lambda}|_{00}=q_1, \quad B_{\lambda}|_{10}=r_1, \quad
B_{\lambda}|_{21}=\sigma
\]
and all the other restrictions are trivial. In particular, the
automorphisms $A$ and $B_{\lambda}$ are finite-state and Sushchansky
group $G_{\lambda}$ is generated by two finite initial automata.
Denote the union of these two automata by $\A_{u,v}$. Its structure
is shown in Figure~\ref{aut_general}. The particular automaton for
$p=3$ and the lexicographic order on $\{(\alpha,\beta) |
\alpha,\beta\in \Fp\}$ is given in Figure~\ref{aut_particular} (all
the arrows not shown in the figures go to the trivial state $1$).

Notice that the word $v$ cannot be periodic since it contains
exactly $p-1$ zeros and $p-1\nmid p^2$. On the contrary $u$ may be
periodic with period $p$. In this case we have $q_i=q_{i+p}$ and the
minimization of $\A_{u,v}$ contains $p^2+2p+5$ states. If $u$ is not
periodic then $\A_{u,v}$ contains $2p^2+p+5$ states. Let $t$ be the
length of the minimal period in $u$ (thus either $t=p$ or $t=p^2$).

%\begin{lemma}
%Properties of $u$ and $v$.
%\begin{enumerate}
%    \item $u_i=1$ or $v_i=1$ for every $i$.
%    \item
%    \item
%\end{enumerate}
%\end{lemma}

\begin{lemma}\label{lemma qr elem abelian group}
The group $\langle q_1,\ldots,q_{t},r_1,\ldots,r_{p^2}\rangle$ is
elementary abelian $p$-group.
\end{lemma}
\begin{proof}
%As we already observed all $q_i,r_j$ are different. Really,
%$q_i\neq q_j$ ($r_i\neq r_j$), since otherwise $u$ ($v$) will be a
%periodic word; $q_i\neq r_j$, because $u\neq v$.

All $q_i,r_j$ have order $p$ because
\[
q_i^p=(q_{i+1}^p,1,1,\ldots,1), \qquad
r_i^p=(r_{i+1}^p,1,1,\ldots,1),
\]
and therefore $q_i^p$ and $r_i^p$ act trivially on the tree.

All $q_i, r_j$ commute with each other, because
\begin{eqnarray*}
q_iq_j=(q_{i+1}q_{j+1},\sigma^{u_i+u_j},1,\ldots,1), \quad
q_jq_i=(q_{j+1}q_{i+1},\sigma^{u_i+u_j},1,\ldots,1);\\
r_ir_j=(r_{i+1}r_{j+1},\sigma^{v_i+v_j},1,\ldots,1), \quad
r_jr_i=(r_{j+1}r_{i+1},\sigma^{v_i+v_j},1,\ldots,1);\\
q_ir_j=(q_{i+1}r_{j+1},\sigma^{u_i+v_j},1,\ldots,1), \quad
r_jq_i=(r_{j+1}q_{i+1},\sigma^{u_i+v_j},1,\ldots,1),
\end{eqnarray*}
so the corresponding pairs act equally on the tree.
\end{proof}

The last lemma implies that the order of $B_{\lambda}$ is $p$. Since
\[
A^{p}=(\sigma^{\frac{p(p-1)}{2}}, \sigma^{\frac{p(p-1)}{2}}, \ldots,
\sigma^{\frac{p(p-1)}{2}})
\]
and $p$ is odd, the order of $A$ is also $p$.

\section{Actions on rooted trees}\label{sec_trees}
Here we describe the structure of the action of $G_\lambda$ on a
$p$-ary tree by means of the orbit tree. This notion is defined
in~\cite{serre:trees} and used in~\cite{gawron_ns:conjugation} to
establish a criterion determining when two automorphisms of a rooted
tree are conjugate. Here we use it to simplify the definition of
Sushchansky groups and show that they admit a faithful
level-transitive action on a regular rooted tree.

\begin{defin}
Let $G$ be a group acting on a regular $p$-ary tree $X^{*}$. The
\emph{orbit tree} of $G$ is a graph whose vertices are the orbits of
$G$ on the levels of $X^{*}$ and two orbits are adjacent if and only
if they contain vertices that are adjacent in $X^{*}$.
\end{defin}

\begin{prop}\label{prop orbit tree}
The structure of the orbit tree of $G_\lambda$ does not depend on
the type $\lambda$ and is shown in Figure~\ref{fig orbit tree}.
\begin{figure}[h]
\begin{center}
\includegraphics{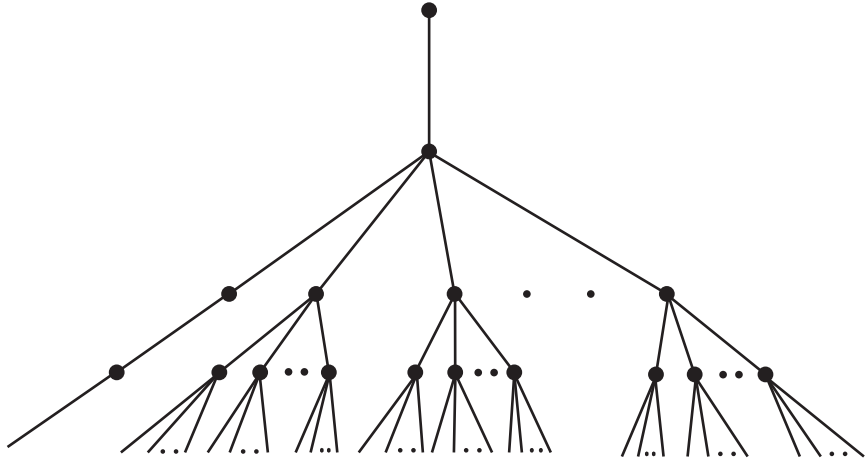}
\caption{The Orbit Tree of Sushchansky group\label{fig orbit tree}}
\end{center}
\end{figure}
\end{prop}

\begin{proof}
Let $T_O$ be the orbit tree of $G_\lambda$. Denote by $\Orb(w)$ the
orbit of the word $w\in X^{*}$ under the action of $G_{\lambda}$.
Define the set
\begin{equation}
V=\left\{ xyw\in X^{*} | xy\in \Orb(00) \mbox{ and } w\in
X^{*}\right\}\cup\{\emptyset\},
\end{equation}
where $\emptyset$ is the root of the tree.

The generator $B_\lambda$ stabilizes the second level of the tree
and hence the orbit $\Orb(00)$ coincides with the orbit of $00$
under the action of the group generated by $A$. The set $V$ and its
compliment $W=X^{*}\setminus V$ are invariant under the action of
$G_{\lambda}$.

Notice that $\{00,10,21\}\subset \Orb(00)$ and the generator
$B_{\lambda}$ acts trivially on all words that lie in the set $W$.
Since the restrictions of $A$ on all words of length $\geq 2$ are
trivial, every element $g\in G_{\lambda}$ that acts trivially on the
second level of the tree must stabilize all the vertices of the set
$W$. Hence, the orbits of $G_\lambda$ on $W$ coincide with the ones
of $A$. Automorphism $A$ acts transitively on the first level and
has order $p$. Therefore the orbit of any word $w\in W$ consists of
$p$ vertices, namely the images of $w$ under the action of the
cyclic group of order $p$ generated by $A$. Therefore the first two
levels of $T_O$ are exactly as shown in Figure~\ref{fig orbit tree}
and $p-1$ vertices on the second level of $T_O$ are the roots of
regular $p$-ary trees.

Let us prove that $G_{\lambda}$ acts transitively on the levels of
the set $V$, i.e. for every $n\geq 1$ the group $G_{\lambda}$ acts
transitively on the set
\[
V_n=\{xyw\in X^{n+1} | xy\in \Orb(00) \mbox{ and } w\in X^{n-1}\}.
\]
We use induction on $n$. For $n=1$ there is nothing to prove. Assume
$G_\lambda$ acts transitively on $V_n$ and consider the $(n+1)$-th
level. Since by construction either $u_{n-1}=1$ or $v_{n-1}=1$, the
restriction of $B_\lambda$ on either $00\ldots01$ or $10\ldots01$ is
equal to $\sigma$. Denote this word as $s$ (here $s\in V_n$) and
note that $B$ stabilizes $s$. To prove the induction step it
suffices for an arbitrary $s'z'\in V_{n+1}$, where $s'\in V_n$ and
$z'\in X$, to construct an element $g\in G_\lambda$ such that
$g(s0)=s'z'$. By the inductive assumption there is an element $h\in
G_\lambda$ such that $h(s)=s'$. Suppose $h^{-1}(s'z')=sz$ for some
letter $z\in X$. Then for $g=(B_\lambda)^{z}h$ (here we consider $z$
as an integer) we have
\begin{eqnarray*}
g(s0)&=&h((B_\lambda)^{z}(s0))=h(s(B_\lambda)^z|_s(0))=h(s(B_\lambda|_s)^z(0))=\\
&=&h(s\sigma^z(0))=h(sz)=s'z'
\end{eqnarray*}
as required.
\end{proof}

The set $V$ has a natural structure of a rooted $p$-ary tree $T$,
where the root $\emptyset$ is connected by an edge with every vertex
in $\Orb(00)$ and there is an edge between $w$ and $wx$ for all
$w\in V$ and $x\in X$. In other words, there is a natural 1-to-1
correspondence between $V$ and vertices of $T$ given by $xyw\mapsto
xw$ for $xy\in\Orb(00)$ and $w\in X^*$. Since the set $V$ is
invariant under the action of $G_{\lambda}$, the group $G_{\lambda}$
acts by automorphisms on the tree $T$. This action has simpler
structure and the following proposition holds.

\begin{prop}\label{prop_simplif}
The action of Sushchansky group $G_{\lambda}$ on the tree $T$ is
faithful, level transitive and has the following form
\begin{equation}
\label{eqn_ABqr}
\begin{array}{rcl}
A &=& \sigma,\\
B_\lambda &=& (q_1,r_1,\sigma,1,\ldots,1),\\
q_i &=& (q_{i+1},\sigma^{u_i},1,\ldots,1),\\
r_i &=& (r_{i+1},\sigma^{v_i},1,\ldots,1).
\end{array}
\end{equation}
\end{prop}
\begin{proof}
The expressions~\eqref{eqn_ABqr} follow directly from the definition
of Sushchansky groups.

Let us prove that this action is faithful. Take an arbitrary
nontrivial element $g\in G_\lambda$. If $g$ acts non-trivially on
the second level of $X^{*}$, then the exponent of $A$ in $g$ is not
divisible by $p$. But then $g$ acts non-trivially on the first level
of $T$ as well because it is fixed under $B_{\lambda}$ and $A$ acts
there as $\sigma$. If $g$ acts trivially on the second level of
$X^{*}$ then it acts trivially on the complement of $V$ in $X^{*}$
according to Proposition~\ref{prop orbit tree}. Therefore to be
nontrivial it must act nontrivially on $T$.

We proved in Proposition~\ref{prop orbit tree} that $G_{\lambda}$
acts transitively on every set $V_n$, which is precisely the $n$th
level of the tree $T$.
\end{proof}

\section{Self-similar closure}
\label{sec_closure}

The Sushchansky group $G_\lambda$ is not generated by all the states
of $\A_{u,v}$ and is not self-similar (see definition below).
However, we can embed it into a larger self-similar group where we
can use some known techniques to derive some important results about
$G_\lambda$ itself. In particular that $G_\lambda$ is amenable
(Corollary~\ref{cor_amenable}) and that the word problem is solvable
in polynomial time (Corollary~\ref{cor_wp}). For the definitions not
given here and more information on self-similar groups we refer
to~\cite{nekrash:self-similar} and~\cite{bar_gs:branch}.

\begin{defin}
A group $G<\Aut X^{*}$ is called \emph{self-similar} if $g|_u\in G$
for any $g\in G$ and word $u\in X^*$. The \emph{self-similar}
closure of $G<\Aut X^*$ is the group generated by all the
restrictions of all the elements of $G$ on words in $X^*$.
\end{defin}

Let $\tilde G_{\lambda}$ be the self-similar closure of
$G_{\lambda}$, i.e. $\tilde G_{\lambda}$ is generated by all the
states of the automaton $\A_{u,v}$. Consider also the self-similar
subgroup $K=\langle
q_1,\ldots,q_{t},r_1,\ldots,r_{p^2},\sigma\rangle$ of $\tilde
G_{\lambda}$.

\begin{lemma}
\label{inf_order}
The group $K$ is not periodic.
\end{lemma}
\begin{proof}
First, consider the case $t=p$. Then all $u_i$'s are equal to $1$
except one equal to $0$. In particular, $\sum_{i=1}^pu_i=p-1$. Then
the element $g=q_1q_2\cdots q_t\sigma^{p-1}$ has representation
\[
g=(q_1q_2\cdots q_t,\sigma^{p-1},1,\ldots,1)\sigma^{p-1}.
\]
Therefore
\[
g^p=(q_1q_2\cdots q_t\sigma^{p-1},*,\ldots,*)=(g,*,\ldots,*).
\]
Since $g$ is nontrivial it must have infinite order.

In case $t=p^2$, exactly $p$ of $u_i$'s are zeros. We mark the
vertices of the cycle of $q_i$'s in the automaton by the
corresponding $u_i$'s. There are at most $p\choose 2$ different
distances between the zeros in the cycle. But the length of the
cycle is $p^2$ so there are
\[\frac{p^2-1}2>\frac{p^2-p}2={p\choose 2}\]
possible distances in the cycle, so let $d$ be a distance that is
not attained as a distance between two zeros.

Now consider the element $g=q_1q_{d+1}\sigma^{u_{p^2}+u_d}$. It can
be written as
\[g=(q_2q_{d+2},\sigma^{u_1+u_{d+1}},1,\ldots,1)\sigma^{u_{p^2}+u_d}.\]
Since the distance between states $q_{p^2}$ and $q_d$ in the cycle
is exactly $d$ at least one of $u_{p^2}$ and $u_d$ is nonzero so
$\sigma^{u_{p^2}+u_d}$ is a cycle of length $p$. Hence
\[g^p=(q_2q_{d+2}\sigma^{u_1+u_{d+1}},*,\ldots,*).\]
Therefore if the order $|g|$ of $g$ is finite, then it is not
smaller than ${p\cdot|q_2q_{d+2}\sigma^{u_1+u_{d+1}}|}$.

Now we repeat this procedure $p^2$ times and on the $i$-th iteration
we get
\[q_{i}q_{d+i}\sigma^{u_{i-1}+u_{d+i-1}}=
(q_{i+1}q_{d+i+1},\sigma^{u_{i}+u_{d+i}},1,
\ldots,1)\sigma^{u_{i-1}+u_{d+i-1}}.\]

Again, the distance between $q_{i-1}$ and $q_{d+i-1}$ is exactly $d$
so $\sigma^{u_{i-1}+u_{d+i-1}}$ is a cycle of length $p$ and
\[(q_{i}q_{d+i}\sigma^{u_{i-1}+u_{d+i-1}})^p=
(q_{i+1}q_{d+i+1}\sigma^{u_{i}+u_{d+i}},*,\ldots,*).\] Therefore
\[|q_{i}q_{d+i}\sigma^{u_{i-1}+u_{d+i-1}}|\geq p\cdot|q_{i+1}q_{d+i+1}\sigma^{u_{i}+u_{d+i}}|.\]
But after $p^2$ steps we will meet $g$ again. So its order cannot be
finite.
\end{proof}

\begin{defin}
A group $G$ acting on the tree $X^*$ is called \emph{weakly regular
branch} over its subgroup $P$, if
\begin{enumerate}
    \item $G$ acts transitively on each level $X^n$, $n \geq 0$;
    \item $P\succ P\times P\times\cdots\times P$ as
    geometric embedding induced by the restriction on some level $X^k$.
\end{enumerate}
In case if $P$ is a subgroup of finite index in $G$, the group $G$
is said to be \emph{regular branch} over $P$.
\end{defin}

\begin{prop}
$\tilde G_{\lambda}$ is a weakly regular branch group over $K^p$.
\end{prop}

\begin{proof}
First of all note that Lemma~\ref{inf_order} guarantees that $K^p$
is nontrivial.  At least one (in fact more) of the $u_i$'s is non
zero, say $u_1$. Then the relations~(\ref{eqn def aut q r}) and
\[
\sigma q_1 \sigma^{p-1}=(\sigma^{u_1},1,\ldots,1,q_2)
\]
show that the set of restrictions of the elements of $K$, that
stabilize the first level $X$ of the tree, on letter $0$ includes
the generators of $K$ and hence the whole group $K$ (therefore
conjugating by $\sigma\in K$ yields that $K$ is self-replicating,
i.e. for any $x\in X$ the projection of $\St_x(K)$ onto the vertex
$x$ coincides with $K$). Thus for any $v\in K$ there is $w\in K$ of
the form
\[w=(v,\sigma^i,1,\ldots,1,q_2^j)\]
for some $i$ and $j$. But then by Lemma~\ref{lemma qr elem abelian
group}
\[
w^p=(v^p,\sigma^{ip},1,\ldots,1,q_2^{jp})=(v^p,1,\ldots,1).
\]
Therefore $K^p\succ K^p\times1\times\cdots\times1$. Since $\sigma$
acts transitively on the first level and belongs to the normalizer
of $K^p$ in $K$ (because
$\sigma^{-1}v^p\sigma=(\sigma^{-1}v\sigma)^p$) by conjugation we get
\[
K^p\succ K^p\times K^p\times\cdots\times K^p,
\]
as geometric embedding.

The transitivity of $\tilde G_\lambda$ on levels follows from the
fact that its subgroup $K$ acts nontrivially on the first level and
is self-replicating, and hence, level transitive. Another
explanation comes from the known fact that a self-similar subgroup
of $\wr_{i\geq1}C_p^{(i)}$ acts level-transitively if and only it is
infinite (see~\cite{bondarenko_gkmnss:clas32}). The proof of the
last fact is similar to the proof of transitivity in
Proposition~\ref{prop orbit tree}.
\end{proof}

We summarize some general properties of $\tilde G_{\lambda}$ in the
following proposition:

\begin{prop}
The self-similar closure of $G_{\lambda}$ is neither torsion, nor
torsion free, level-transitive group of tree automorphisms.
Moreover, it is generated by a bounded automaton, hence it is
contracting and amenable.
\end{prop}
\begin{proof}
The first three assertions are already proved above. The automaton
$\A_{u,v}$ is bounded by Corollary 14 in \cite{sid:cycl} (see the
definition there as well). As a corollary $\tilde G_{\lambda}$ is
contracting (see \cite{bond_n:pcf}) and amenable (see
\cite{bknv:amenab}).
\end{proof}
\begin{cor}
\label{cor_amenable}
$G_{\lambda}$ is amenable.
\end{cor}
Note also that the last corollary follows from Theorem~\ref{thm
SuschGrowth}.

\begin{cor}
\label{cor_wp}
The word problem in $G_{\lambda}$ is solvable in polynomial time.
\end{cor}
\begin{proof}
See Proposition 2.13.10 in \cite{nekrash:self-similar}.
\end{proof}

\section{Intermediate growth}
\label{sec_growth}

Let $G$ be a group finitely generated by a set $S$. \textit{The
growth function of $G$} is defined by
\[
\gamma_G(n)=\left|\{g\in G |  g=s_1s_2\ldots s_k \mbox{ for some }
s_i\in S\cup S^{-1}, k\leq n\}\right|.
\]
Two functions $\gamma_1$ and $\gamma_2$ are called
\textit{equivalent} if there exists a constant $C>0$ such that
$\gamma_1(\frac{1}{C}n)\leq \gamma_2(n)\leq\gamma_1(Cn)$ for all
$n$. The growth function $\gamma_G$ depends both on $G$ and on $S$,
but the equivalence class of $\gamma_G$ does not depend on $S$.

In 1968 John Milnor asked about the existence of finitely generated
groups with growth that is intermediate between polynomial and
exponential. The first examples of such groups were provided by
R.I.~Grigorchuk in~\cite{grigorch:milnor}, where he constructed
uncountable family of such groups. In particular, it was shown, that
there are groups of intermediate growth generated by automata with
$5$ states, namely, $G_\omega$ for $\omega=(012)^\infty$ (not to be
confused with Sushchansky groups $G_\lambda$). These examples were
generalized to the notion of $\G$ groups \cite{bar_gs:branch}. Under
some finiteness restriction all $\G$ groups have intermediate
growth.

Recently it was proved~\cite{bux_p:iter_monodromy} that there is a
$4$-state automaton over a $2$-letter alphabet generating a group of
intermediate growth. This group itself is isomorphic to the iterated
monodromy group of the map $f(z)=z^2+i$. But it is still an open
question whether there is a group of intermediate growth generated
by a $3$-state automaton over a $2$-letter alphabet.

In view of the examples above it is not very surprising that the two
of the pioneering examples of infinite finitely generated periodic
groups introduced by S.V.~Aleshin in~\cite{aleshin:burnside} and
V.I.~Sushchansky in~\cite{sushch:burnside} also have intermediate
growth. For Aleshin group it follows from the intermediate growth of
Grigorchuk group and the result of
Y.I.~Merzlyakov~\cite{merzlyakov:periodic}, who proved that Aleshin
group contains a subgroup of finite index isomorphic to the
subdirect product of four copies of Grigorchuk group. Also the
relation between these two groups was studied
in~\cite{grigorch:habil}.

As was mentioned above in~\cite{grigorch:degrees85} R.I.~Grigorchuk
pointed out that all Sushchansky groups have intermediate growth,
but only the idea of proof was given. In this paper we give a
complete proof of this fact based on the results
from~\cite{bartholdi_s:growth}.

At the present moment the main method of obtaining the upper bounds
for growth functions of groups was originated by R.I.~Grigorchuk
in~\cite{grigorch:degrees}. Different modifications of this method
in~\cite{bartholdi:growth,muchnik_p:growth,bartholdi_s:growth}
allowed to improve existing estimates and to prove the estimates for
new groups.

As for the lower bounds for growth functions, there are several
techniques. In~\cite{grigorch:degrees} R.I.~Grigorchuk uses
self-similarity to obtain the lower bound of the form $e^{\sqrt n}$
for most of his groups. Moreover, he shows that any group $G$ that
is abstractly commensurable with its own power $G^k$ for some
$k\geq2$ has a growth function not smaller that $e^{n^\alpha}$ for
some $0<\alpha\leq1$.

In~\cite{grigorch:hilbert} R.I.~Grigorchuk used bounds on the
coefficients of Hilbert-Poincar\'e series of graded algebras
associated with groups to bound their growth functions. Namely, it
was obtained that any residually $p$-group whose growth function is
not bounded above by polynomial, must grow at least as $e^{\sqrt
n}$.

Y.G.~Leonov~\cite{leonov:bound}, L.~Bartholdi and
Z.~\v{S}uni\'c~\cite{bartholdi:growth, bartholdi_s:growth} used more
advanced techniques (common in spirit to the ones used
in~\cite{grigorch:degrees}) also based on certain self-similarity of
the groups acting on trees. In obtaining the lower bounds for the
growth functions of these groups the important role was played by
the property, which is in some sense opposite to contraction. The
main idea is that the restrictions of elements can not be much
shorter than the elements themselves.

A.~Erschler used random walks and Poisson boundary to approach to
this question. In particular, in~\cite{erschler:subexp} it was shown
that the growth function of Grigorchuk group $G_\omega$ for
$\omega=(01)^\infty$, which is generated by 5-state automaton, grows
faster than $e^{n^\alpha}$ for any $\alpha<1$. The upper estimate of
the same sort was obtained for this group in spirit
of~\cite{grigorch:degrees}, which shows that groups $G_\omega$ for
$\omega=(012)^\infty$ and $\omega=(01)^\infty$ have essentially
different growth functions.

Recall the definition of a $\G$ group.

\begin{defin}
Let $R$ be a subgroup of $Sym(X)$, $D$ be any group with a sequence
of homomorphisms $w_i:D\to\Sym(X)$, $i\geq 1$. Then $R$ acts on the
first level of $X^{*}$ and $D$ acts on $X^{*}$ in the following way.
Each $d\in D$ defines the automorphism $\hat{d}$ that acts trivially
on the first level and is given by its restrictions
\[
\hat{d}\bigl|_{0^i1}=w_i(d), i\geq1
\]
and all the other restrictions act trivially on $X$. Denote
$\hat{D}=\{\hat{d} \mid d\in D\}$.

The group $G=\langle R, \hat{D}\rangle$ is called a \emph{$\G$
group} if the following conditions are satisfied:
\begin{enumerate}
\item[(i)] The groups $R$ and $w_i(D), i\geq1$, act transitively on
$X$.
\item[(ii)] For each $d\in D$ the permutation $w_i(d)$ is trivial for infinitely many
indices.
\item[(iii)] For each nontrivial $d\in D$ the permutation $w_i(d)$ is nontrivial for infinitely many
indices.
\end{enumerate}
The groups $R$ and $D$ are called \textit{the root part} and
\textit{the directed part} of $G$ correspondingly.
\end{defin}

Note that in~\cite{bar_gs:branch} the definition of a $\G$ group is
given in slightly more general settings. The results
in~\cite{bartholdi_s:growth} and~\cite{bar_gs:branch} imply the
following theorem.

\begin{thm}\label{thm growth G groups}
All $\G$ groups with finite directed part have intermediate growth.
\end{thm}

There is a lower bound for the growth of such groups given
in~\cite{bar_gs:branch}:
\begin{equation}
\label{eqn lower bound growth G groups} \gamma_G(n) \succeq
e^{n^\alpha},
\end{equation}
where $\alpha=\frac{\log(|X|)}{\log(|X|)+\log(2)}$.

The sequence of homomorphisms $w_i$ in the definition of a $\mathsf
G$ group is called $r$-homogeneous, if for every finite subsequence
of $r$ consecutive homomorphisms $w_{i+1},w_{i+2},\ldots,w_{i+r}$
every element of $D$ is sent to the identity by at least one of the
homomorphisms from this finite subsequence. In particular, if the
sequence of homomorphisms $\{w_i,i\geq1\}$ defining a $\mathsf G$
group is periodic with period $r$, it is also $r$-homogeneous.

It is proved in~\cite{bartholdi_s:growth} that in case of
$r$-homogeneous sequence of defining homomorphisms there is an
estimate of the upper bound on the growth function. Moreover, in
this case if the directed part has finite exponent there is an upper
bound on the torsion growth function $\pi(n)$ (the maximal order of
an element of length at most $n$).

\begin{thm}[$\eta$-estimate]
\label{thm_eta-estimate} Let $G$ be a $\mathsf G$ group defined by
an $r$-homogeneous sequence of homomorphisms. Then the growth
function of the group $G$ satisfies
\begin{equation}
\label{eqn upper bound growth G groups} \gamma_G(n) \preceq
e^{n^\beta},
\end{equation}
where $\beta=\frac{\log(|X|)}{\log(|X|)-\log(\eta_r)}<1$ and
$\eta_r$ is the positive root of the polynomial
$x^r+x^{r-1}+x^{r-2}-2$.

If the directed part $D$ of $G$ has finite exponent $q$, then the
group $G$ is torsion and there exists a constant $C>0$, such that
the torsion growth function satisfies
\begin{equation}\label{eqn_upper_bound_period_growth}
\pi(n)\leq Cn^{\log_{1/\eta_r}(q)}.
\end{equation}

%Also if the directed part $D$ is of finite exponent $q$, then the
%torsion growth function $\pi_G(n)$ (maximal order of an element of
%length up to $n$) satisfies
%\[\pi_G(n)\leq Cn^{\log_{1/\eta_r}q}\]
%for some positive constant $C$.
\end{thm}

Sushchansky groups $G_{\lambda}$ are not $\G$ groups, because the
automorphism $B_{\lambda}$ cannot be expressed as $\hat{d}$ for some
homomorphisms $w_i$. On the other hand, the automorphisms $q_i$ and
$r_i$ can, and the following proposition shows that the self-similar
closure of $G_{\lambda}$ contains a subgroup which is a $\G$ group.
Since the simplified definition of $G_\lambda$ from
Proposition~\ref{prop_simplif} does not simplify considerably the
proofs in this section, we use the original definition in order to
make this section independent.

\begin{prop}\label{prop H is G group}
The group $H=\langle q_1, r_1,\sigma \rangle$ is a $\G$ group with
finite directed part defined by a periodic sequence of homomorphisms
with period $p^2$.
\end{prop}

\begin{proof}
We prove that the subgroups $\langle q_1,r_1\rangle$ and
$\langle\sigma\rangle$ are the directed and the root parts of $H$.

First observe that $\langle q_1,r_1\rangle\simeq\Z_p\oplus\Z_p$.
Indeed, the group $\langle q_1,r_1\rangle$ is elementary abelian
$p$-group by Lemma \ref{lemma qr elem abelian group}. Suppose that
$r_1\in\langle q_1\rangle$, $r_1=q_1^k$. Comparing restrictions on
words $0\ldots 01$ we get $v_i=ku_i$. Contradiction, since $u_i=0$
and $v_i=1$ for $i$ with $\beta_i=0$.

Consider the periodic sequence of homomorphisms $w_i:\langle
q_1,r_1\rangle\to\Sym(X)$ with period $p^2$ given by
$w_i(q_1)=\sigma^{u_i}$ and $w_i(r_1)=\sigma^{v_i}$. Then for any
$d\in\langle q_1,r_1\rangle$ the associated $\hat{d}$ from the
definition of a $\G$ group coincides with the automorphism $d$. To
complete the proof we need to check the conditions (i)--(iii) from
the definition of a $\G$ group.

(i) The root part generated by $\sigma$ acts transitively on $X$.
Furthermore, for any $i\geq 1$
\begin{eqnarray*}
w_i(q_1)=\sigma, &\mbox{ if } \beta_i\neq 0;\\
w_i(r_1)=\sigma, &\mbox{ if } \beta_i= 0.
\end{eqnarray*}
In any case $w_i(\langle q_1,r_1\rangle)$ contains $\sigma$ and thus
acts transitively on $X$.

(ii),(iii) Let $d=q_1^kr_1^l$, $k,l\in\Z_p$, be an arbitrary
nontrivial element of $\langle q_1,r_1\rangle$. Since the sequence
$w_i$ is periodic it suffices to show at least one occurrence of
trivial and one occurrence of nontrivial $w_i(d)$.

Find $i$ such that
\begin{eqnarray*}
(\alpha_i,\beta_i)=(1,0), &\mbox{ if } l=0;\\
(\alpha_i,\beta_i)=(k,l), &\mbox{ if } l\neq 0.
\end{eqnarray*}
Then
\[
w_i(d)=\left\{%
\begin{array}{ll}
    w_i(q_1^k)=\sigma^{ku_i}=1, & \hbox{ if } l=0; \\
    w_i(q_1^kr_1^l)=\sigma^{ku_i+lv_i}=\sigma^{k+l(-k/l)}=1, & \hbox{ if } l\neq 0. \\
\end{array}%
\right.
\]

For a nontrivial occurrence find $i$ such that
\begin{eqnarray*}
(\alpha_i,\beta_i)=(0,1), &\mbox{ if } l=0;\\
(\alpha_i,\beta_i)=(1,0), &\mbox{ if } l\neq 0.
\end{eqnarray*}
Then
\[
w_i(d)=\left\{%
\begin{array}{ll}
    w_i(q_1^k)=\sigma^{ku_i}=\sigma^k, & \hbox{ if } l=0; \\
    w_i(q_1^kr_1^l)=\sigma^{ku_i+lv_i}=\sigma^l, & \hbox{ if } l\neq 0. \\
\end{array}%
\right.
\]
\end{proof}

The last proposition shows that the growth function of $H$ satisfies
inequalities~(\ref{eqn lower bound growth G groups}) and~(\ref{eqn
upper bound growth G groups}), for $r=p^2$. Also note that it is
proved in~\cite{bar_gs:branch} that a $\G$ group is torsion if and
only if its directed part $D$ is torsion. Therefore, the group $H$
is torsion. The next proposition exhibits another branch structure
inside $\tilde G_\lambda$.

\begin{prop}\label{prop H is branch}
The group $H=\langle q_1, r_1,\sigma \rangle$ is regular branch over
its commutator subgroup $H'$.
\end{prop}

\begin{proof}
Let $H_k=\langle q_k, r_k,\sigma \rangle$, $k=1,\ldots,p^2$ be the
subgroups of $\tilde G_\lambda$. First we show that
\begin{equation}
\label{eqn:H_k_branch} H_k'\succeq H_{k+1}'\times H_{k+1}'\times\dots\times H_{k+1}'
\end{equation}
for all $k$. Indeed, at least one of $u_k$ and $v_k$ is nonzero.
Suppose $u_k\neq 0$. Then relations
$q_k=(q_{k+1},\sigma^{u_k},1,\ldots,1)$ and
$r_k=(r_{k+1},\sigma^{v_k},1,\ldots,1)$ imply
\begin{eqnarray*}
~[q_k,r_k]&=&([q_{k+1},r_{k+1}],1,\ldots,1),\\
~[q_k,(q_k^{\sigma^{-1}})^{1/u_k}]&=&([q_{k+1},\sigma],1,\ldots,1),\\
~[r_k,(q_k^{\sigma^{-1}})^{1/u_k}]&=&([r_{k+1},\sigma],1,\ldots,1).
\end{eqnarray*}
Since the projection of the stabilizer of the first level in $H_k$
on the leftmost vertex coincides with $H_{k+1}$ we get $H_k'\succeq
H_{k+1}'\times 1\times\dots\times 1$. Conjugation by $\sigma\in H_k$
implies inclusion~\eqref{eqn:H_k_branch}. Since $H_1=H_{p^2+1}=H$,
we obtain $H'\succeq H'\times H'\times\dots\times H'$ as geometric
embedding induced by the restriction on $X^{p^2}$.

The transitivity of $H$ on the levels is proved by the method used
in Proposition~\ref{prop orbit tree}.

Now $H$ is a torsion $p$-group, hence, so is $H/H'$, which is
abelian. But each torsion finitely generated abelian group is
finite. Thus, $H'$ is a subgroup of finite index in $H$.
\end{proof}

When we deal with a group $G$ of automorphisms of $X^{*}$, it is
sometimes difficult to say something about the whole group, but we
know something about the group $P$ generated by all the restrictions
of the elements in $G$ on some level $k$ of the tree. In case $G$ is
self-similar, $P$ is a subgroup of $G$ and if $G$ is
self-replicating, $P$ coincides with $G$. Some properties of $P$ are
inherited by $G$ itself. In particular, if $P$ is finite or torsion
then so is $G$ (the converse is not true). But what we are
interested in here is that the growth of $G$ can be estimated in
terms of the growth of $P$.

Let $S$ be a finite generating set of $G$. Then  $P$ is generated by
the set $\tilde{S}$ of the restrictions of all elements of $S$ on
all vertices of $k$-th level $X^k$ of the tree. The following lemma
holds.

\begin{lemma}\label{lemma growth estimate}
The growth function $\gamma_G(n)$ of the group $G$ with respect to
$S$ is bounded from above by
\begin{equation}\label{eqn in lemma growth estimate}
\gamma_G(n)\preceq\bigl(\gamma_P(n)\bigr)^{|X|^k},
%(p!)^{\frac{1-p^k}{1-p}}\bigl(\gamma_P(n)\bigr)^{p^k},
\end{equation}
where $\gamma_P(n)$ is the growth function of the group $P$ with
respect to $\tilde{S}$. In particular, the growth type of $G$
(finite, polynomial, intermediate or exponential) cannot exceed the
one of $P$.
\end{lemma}

\begin{proof}
Let $g\in G$ be an element of length $n$ with respect to the
generating set $S$. This element induces a permutation $\pi_k$ of
the $k$-th level of the tree and $|X|^k$ restrictions $g|_{v}, v\in
X^k$, on words of length $k$. Moreover, different automorphisms
correspond to different tuples $(\pi_k,\{g|_v, v\in X^k\})$ of
restrictions and permutations. Each such a restriction is a word of
length not greater than $n$ with respect to the generating set
$\tilde{S}$ of $P$. So for each vertex $v\in X^k$ the number of
possible restrictions on $v$ is bounded from above by $\gamma_P(n)$.
\end{proof}

The following corollary shows an easy way to construct new examples
of groups with intermediate (finite, polynomial, exponential)
growth.

\begin{cor}
Let $F$ be a finite set of automorphisms from $\Aut X^{*}$, whose
restrictions on some level $k$ belong to $G$ (in particular, $F$
could be a set of finitary automorphisms). Then
\[
\gamma_G(n)\precsim\gamma_{\langle
G,F\rangle}(n)\precsim\bigl(\gamma_G(n)\bigr)^{|X|^k}.
\]
where $\gamma_{\langle G,F\rangle}(n)$ is the growth function of the
group $\langle G,F\rangle$ with respect to $S\cup F$.
\end{cor}

In particular the previous corollary shows that if a group $G$ is
generated by a finite automaton, then the growth type of this group
depends only on the nucleus (see definition
in~\cite{nekrash:self-similar}) of this automaton.

An interesting question is whether it is true that if $G$ grows
faster than polynomially then $\gamma_G(n)\thicksim\gamma_{\langle
G,F\rangle}(n)$.

We are ready to prove the main results.

\begin{thm}\label{thm SuschGrowth}
All Sushchansky $p$-groups have intermediate growth. The growth
function of each Sushchansky $p$-group $G_\lambda$ satisfies
\[e^{n^\alpha}\preceq\gamma_{G_\lambda}(n)\preceq e^{n^\beta},\]
where $\alpha=\frac{\log(p)}{\log(p)+\log(2)}$,
$\beta=\frac{\log(p)}{\log(p)-\log(\eta_r)}$ and $\eta_r$ is the
positive root of the polynomial $x^r+x^{r-1}+x^{r-2}-2$, where
$r=p^2$.
\end{thm}
\begin{proof}
%Now it is enough to consider the projections of $A$ and $B$ on the second level of the tree. We
%have:
%\[A=(1,1,\ldots,1)\pi_A,\]
%\[B=(\underbrace{d,1,\ldots,1}_{p}, \underbrace{f,1,\ldots,1}_{p}, \underbrace{1,\sigma,1,\ldots,1}_{p}, 1,\ldots,1).\]
The group generated by all the restrictions of elements of
$G_{\lambda}$ on the second level is $H=\langle
q_1,r_1,\sigma\rangle$, which is a $\G$ group of intermediate growth
by Proposition~\ref{prop H is G group} and Theorems~\ref{thm growth
G groups} and~\ref{thm_eta-estimate}, whose growth function
satisfies inequalities~(\ref{eqn lower bound growth G groups})
and~(\ref{eqn upper bound growth G groups}). Therefore by
Lemma~\ref{lemma growth estimate} the Sushchansky group
$G_{\lambda}$ has subexponential growth function, which satisfies
inequality
\begin{equation}\label{eqn in thmSuschGrowth est above}
\gamma_G(n)\precsim(\gamma_H(n))^{p^2}\precsim\gamma_H(n).
\end{equation}
The last part of this inequality follows from Proposition~\ref{prop
H is branch}, where it is proved that $H$ is regular branch over
$H'$.

Now consider the subgroup $L=\langle B_{\lambda},
AB_{\lambda}A^{p-1}, A^2B_{\lambda}A^{p-2}\rangle$ of $G_{\lambda}$.
This subgroup stabilizes the second level of the tree and the
restrictions of the generators on the second level look like:
\[
\begin{array}{lcl}
B_{\lambda}&=&(q_1,\ast,\ldots,\ast),\\
AB_{\lambda}A^{p-1}&=&(r_1,\ast,\ldots,\ast),\\
A^2B_{\lambda}A^{p-2}&=&(\sigma,\ast,\ldots,\ast).\\
\end{array}
\]

Each word of length $n$ in $L$ will be projected on the
corresponding word of length $n$ in $H$. Therefore
$\gamma_L(n)\geq\gamma_H(n)$ for all $n\geq1$. But $L$ is a finitely
generated subgroup of $G_\lambda$. Thus
\begin{equation}\label{eqn in thmSuschGrowth est beneath}
\gamma_H(n)\precsim\gamma_L(n)\precsim\gamma_G(n).
\end{equation}
Inequalities~(\ref{eqn in thmSuschGrowth est above}) and~(\ref{eqn
in thmSuschGrowth est beneath}) imply
\begin{equation}
\gamma_G(n)\thicksim\gamma_H(n).
\end{equation}
\end{proof}

Finally, it was mentioned above that the group $H$ is torsion as a
$\G$ group with torsion directed part. But periodicity of $H$
implies that $G_\lambda$ is periodic as well. This gives a different
proof of Theorem~\ref{thm_susch} proved by V.I.~Sushchansky. The
theory of $\G$ groups allows to sharpen this result.

\begin{thm}
There is a constant $C>0$, such that the torsion growth function of
each Sushchansky $p$-group $G_\lambda$ satisfies inequality
\[\pi_{G_\lambda}(n)\leq Cn^{\log_{1/\eta_r}(p)},\]
where $\eta_r$ is the same as in the previous theorem.
\end{thm}

\begin{proof}
By Proposition~\ref{prop H is G group} the group  $H$ is a $\G$
group defined by a $p^2$-homogenous sequence of homomorphisms, whose
directed part $\langle q_1,r_1\rangle$ is an elementary abelian
$p$-group (see Lemma~\ref{lemma qr elem abelian group}). Therefore
by Theorem~\ref{thm_eta-estimate} the torsion growth function
$\pi_H(n)$ satisfies inequality
\[\pi_{H}(n)\leq C_1n^{\log_{1/\eta_r}(p)}\]
for some constant $C_1$.

For any element $g$ of length $n$ in $G_\lambda$, $g^p$ stabilizes
the second level of the tree and the restrictions of $g^p$ at the
vertices of the second level are the elements of $H$, whose length
is not bigger than $pn$. Hence, the order of $g^p$ cannot be bigger
than the least common multiple of the orders of $g|_v$, $v\in X^2$.
Since the orders of these restrictions are the powers of $p$, the
least common multiple coincides with the maximal order among the
restrictions. This implies
\[\mbox{Order}(g)=p\cdot\mbox{Order}(g^p)\leq p\pi_H(pn)\leq
pC_1(pn)^{\log_{1/\eta_r}(p)}\leq Cn^{\log_{1/\eta_r}(p)}\] for
$C=C_1p^{\log_{1/\eta_r}(p)+1}$.

\end{proof}

%\bibliographystyle{alpha}
%\bibliography{../../mylib}

\begin{thebibliography}{BGK{\etalchar{+}}06}

\bibitem[Ale72]{aleshin:burnside}
S.~V. Ale{\v{s}}in.
\newblock Finite automata and the {B}urnside problem for periodic groups.
\newblock {\em Mat. Zametki}, 11:319--328, 1972.

\bibitem[Bar98]{bartholdi:growth}
Laurent Bartholdi.
\newblock The growth of {G}rigorchuk's torsion group.
\newblock {\em Internat. Math. Res. Notices}, (20):1049--1054, 1998.

\bibitem[BGK{\etalchar{+}}06]{bondarenko_gkmnss:clas32}
Ievgen Bondarenko, Rostislav Grigorchuk, Rostyslav Kravchenko,
Yevgen Muntyan,
  Volodymyr Nekrashevych, Dmytro Savchuk, and Zoran \v{S}uni\'{c}.
\newblock Groups generated by $3$-state automata over $2$-letter alphabet, {I}.
\newblock (available at \emph{http://arxiv.org/abs/math.GR/0612178}), 2006.

\bibitem[BG{\v{S}}03]{bar_gs:branch}
Laurent Bartholdi, Rostislav~I. Grigorchuk, and Zoran
{\v{S}}uni{\'k}.
\newblock Branch groups.
\newblock In {\em Handbook of algebra, Vol. 3}, pages 989--1112. North-Holland,
  Amsterdam, 2003.

\bibitem[BKNV06]{bknv:amenab}
Laurent Bartholdi, Vadim Kaimanovich, Volodymyr Nekrashevych, and
Balint Virag.
\newblock Amenability of automata groups.
\newblock (preprint), 2006.

\bibitem[BN03]{bond_n:pcf}
E.~Bondarenko and V.~Nekrashevych.
\newblock Post-critically finite self-similar groups.
\newblock {\em Algebra Discrete Math.}, (4):21--32, 2003.

\bibitem[BP06]{bux_p:iter_monodromy}
Kai-Uwe Bux and Rodrigo P{\'e}rez.
\newblock On the growth of iterated monodromy groups.
\newblock In {\em Topological and asymptotic aspects of group theory}, volume
  394 of {\em Contemp. Math.}, pages 61--76. Amer. Math. Soc., Providence, RI,
  2006.
\newblock (available at \emph{http://www.arxiv.org/abs/math.GR/0405456}).

\bibitem[B{\v{S}}01]{bartholdi_s:growth}
Laurent Bartholdi and Zoran {\v{S}}uni{\'k}.
\newblock On the word and period growth of some groups of tree automorphisms.
\newblock {\em Comm. Algebra}, 29(11):4923--4964, 2001.

\bibitem[Ers04]{erschler:subexp}
Anna Erschler.
\newblock Boundary behavior for groups of subexponential growth.
\newblock {\em Ann. of Math. (2)}, 160(3):1183--1210, 2004.

\bibitem[Glu61]{glushkov:automata}
V.~M. Glu{\v{s}}kov.
\newblock Abstract theory of automata.
\newblock {\em Uspehi Mat. Nauk}, 16(5 (101)):3--62, 1961.

\bibitem[GNS00]{gns00:automata}
R.~I. Grigorchuk, V.~V. Nekrashevich, and V.~I. Sushchanski{\u\i}.
\newblock Automata, dynamical systems, and groups.
\newblock {\em Tr. Mat. Inst. Steklova}, 231(Din. Sist., Avtom. i Beskon.
  Gruppy):134--214, 2000.

\bibitem[GNS01]{gawron_ns:conjugation}
Piotr~W. Gawron, Volodymyr~V. Nekrashevych, and Vitaly~I.
Sushchansky.
\newblock Conjugation in tree automorphism groups.
\newblock {\em Internat. J. Algebra Comput.}, 11(5):529--547, 2001.

\bibitem[Gol64]{golod:p-groups}
E.~S. Golod.
\newblock On nil-algebras and finitely approximable {$p$}-groups.
\newblock {\em Izv. Akad. Nauk SSSR Ser. Mat.}, 28:273--276, 1964.

\bibitem[Gri80]{grigorch:burnside}
R.~I. Grigor{\v{c}}uk.
\newblock On {B}urnside's problem on periodic groups.
\newblock {\em Funktsional. Anal. i Prilozhen.}, 14(1):53--54, 1980.

\bibitem[Gri83]{grigorch:milnor}
R.~I. Grigorchuk.
\newblock On the {M}ilnor problem of group growth.
\newblock {\em Dokl. Akad. Nauk SSSR}, 271(1):30--33, 1983.

\bibitem[Gri84]{grigorch:degrees}
R.~I. Grigorchuk.
\newblock Degrees of growth of finitely generated groups and the theory of
  invariant means.
\newblock {\em Izv. Akad. Nauk SSSR Ser. Mat.}, 48(5):939--985, 1984.

\bibitem[Gri85a]{grigorch:degrees85}
R.~I. Grigorchuk.
\newblock Degrees of growth of {$p$}-groups and torsion-free groups.
\newblock {\em Mat. Sb. (N.S.)}, 126(168)(2):194--214, 286, 1985.

\bibitem[Gri85b]{grigorch:habil}
R.I. Grigorchuk.
\newblock {\em Groups with intermediate growth function and their
  applications}.
\newblock Habilitation, Steklov Institute of Mathematics, 1985.

\bibitem[Gri89]{grigorch:hilbert}
R.~I. Grigorchuk.
\newblock On the {H}ilbert-{P}oincar\'e series of graded algebras that are
  associated with groups.
\newblock {\em Mat. Sb.}, 180(2):207--225, 304, 1989.

\bibitem[Kal48]{kalou:la_structure}
L{\'e}o Kaloujnine.
\newblock La structure des {$p$}-groupes de {S}ylow des groupes sym\'etriques
  finis.
\newblock {\em Ann. Sci. \'Ecole Norm. Sup. (3)}, 65:239--276, 1948.

\bibitem[Leo01]{leonov:bound}
Yu.~G. Leonov.
\newblock On a lower bound for the growth of a 3-generator 2-group.
\newblock {\em Mat. Sb.}, 192(11):77--92, 2001.

\bibitem[Mer83]{merzlyakov:periodic}
Yu.~I. Merzlyakov.
\newblock Infinite finitely generated periodic groups.
\newblock {\em Dokl. Akad. Nauk SSSR}, 268(4):803--805, 1983.

\bibitem[Mil68]{milnor:problem}
J.~Milnor.
\newblock Problem $5603$.
\newblock {\em Amer. Math. Monthly}, 75:685--686, 1968.

\bibitem[MP01]{muchnik_p:growth}
Roman Muchnik and Igor Pak.
\newblock On growth of {G}rigorchuk groups.
\newblock {\em Internat. J. Algebra Comput.}, 11(1):1--17, 2001.

\bibitem[Nek05]{nekrash:self-similar}
Volodymyr Nekrashevych.
\newblock {\em Self-similar groups}, volume 117 of {\em Mathematical Surveys
  and Monographs}.
\newblock American Mathematical Society, Providence, RI, 2005.

\bibitem[Ser03]{serre:trees}
Jean-Pierre Serre.
\newblock {\em Trees}.
\newblock Springer Monographs in Mathematics. Springer-Verlag, Berlin, 2003.

\bibitem[Sid00]{sid:cycl}
S.~Sidki.
\newblock Automorphisms of one-rooted trees: growth, circuit structure and
  acyclicity.
\newblock {\em J. of Mathematical Sciences (New York)}, 100(1):1925--1943,
  2000.

\bibitem[Sus79]{sushch:burnside}
V.~I. Sushchansky.
\newblock Periodic permutation $p$-groups and the unrestricted {Burnside}
  problem.
\newblock {\em DAN SSSR.}, 247(3):557--562, 1979.
\newblock (in Russian).

\end{thebibliography}

\newcommand{\etalchar}[1]{$^{#1}$}
\def\cprime{$'$} \def\cprime{$'$} \def\cprime{$'$} \def\cprime{$'$}

\end{document}